\documentclass[runningheads,a4paper]{llncs}

\usepackage{amsfonts,amssymb,amsmath}
\usepackage{natbib}

\usepackage{wrapfig}
\usepackage{llncsdoc}
\usepackage{epsfig,psfrag}
\usepackage{latexsym}
\usepackage{longtable}
\usepackage{tikz,pgf}
\usepackage{subfig}
\usepackage{wrapfig}

\usepackage{graphicx}
\usepackage{epstopdf}
\usepackage{color}

\newcommand{\bqn}{\begin{eqnarray}}
\newcommand{\eqn}{\end{eqnarray}}
\newcommand{\bq}{\begin{eqnarray*}}
\newcommand{\eq}{\end{eqnarray*}}
\usepackage[ruled,vlined]{algorithm2e}

\usepackage{url}
\usepackage{hyperref}
\hypersetup{colorlinks,%
citecolor=black,%
filecolor=blue,%
linkcolor=red,%
urlcolor=blue,%
pdftex}

\newcommand{\blue}[1]{{\color{blue} #1}}

\begin{document}

\title{Introduction to Random fields}
 \titlerunning{Introduction to Random fields}
 
\author{Moo K. Chung
 }
\institute{
University of Wisconsin-Madison, USA\\
\vspace{0.3cm}
\blue{\tt mkchung@wisc.edu}
}
\authorrunning{Chung}

\maketitle

\begin{center}
March 1, 2007
\end{center}

\pagenumbering{arabic}

\index{random field theory}
\index{fields ! random}

General linear models (GLM) are often constructed and used in statistical inference  at the voxel level  in brain imaging. 
In this paper, we explore the basics of random fields and the multiple comparisons on the random fields, which are
 necessary to properly threshold statistical maps for the whole image at specific statistical significance level. The multiple comparisons are crucial in determining overall statistical significance in correlated test statistics over the whole brain. In practice, $t$- or $F$-statistics in adjacent voxels are correlated. So there is the problem of multiple comparisons, which we have simply neglected up to now. For multiple comparisons that account for spatially correlated test statistics, various methods were proposed: Bonferroni correction, random field theory \citep{ worsley.1994, worsley.1996}, false discovery rates \citep{benjamini.1995, benjamini.2001, genovese.2002} and permutation tests \citep{nichols.2002}. Among them, we will explore the random field approach.

\section{Introduction}

Suppose we measure temperature $Y$ at position $x$ and time $t$ in a classroom $M \in
\mathbb{R}^3$. Since every measurement will be error-prone, we
model the temperature as
$$Y(x,t) = \mu(x,t) + \epsilon(x,t)$$
where $\mu$ is the unknown signal and $\epsilon$ is the
measurement error. The measurement error can be modelled as a
random variable. So at each point $(x,t) \in M \otimes
\mathbb{R}^+$, measurement error $\epsilon(x,t)$ is a random
variable. The collection of random variables $$\{\epsilon(x,t):
(x,t) \in M \otimes \mathbb{R}^+\}$$ is called a {\em stochastic
process}. The generalization of a continuous stochastic process defined in $\mathbb{R}$ to a higher dimensional abstract space indexed by a spatial variable  is called a {\em random field}. For an introduction to random fields, see \citet{adler.2007}, \citet{dougherty.1999} and \citet{yaglom.1987}. A formal measure theoretic definition can be found in \citet{adler.1981} and \citet{gikhman.1996}.

In brain imaging studies, it is necessary to model measurements at each voxel as a random field. For instance, in the deformation-based morphometry (DBM), deformation fields are usually modeled as continuous random fields \citep{chung.2001.ni}. In the random field theory as used in \citep{worsley.1994, worsley.1996}, measurement $Y$ at voxel position $x \in \mathcal{M}$ is modeled  as
$$Y(x) = \mu(x) + \epsilon(x),$$
where $\mu$ is the unknown signal to be estimated and $\epsilon$ is the
measurement error. The measurement error at each fixed $x$ can be modeled as a
random variable. Then the collection of random variables $\{\epsilon(x): x \in \mathcal{M} \}$ is called a {\em stochastic
process} or {\em random field}. The more precise measure-theoretic definition can be found in \citep{adler.2007}. Random field modeling can be done beyond the usual Euclidean space to curved cortical and subcortical manifolds \citep{joshi.1998, chung.2003.ni}.

\section{Random Fields} 
We start with defining a random field more formally using random variables.  Our construction follows from \citet{adler.1981}.

\begin{definition}
Given a probability space, a random field $T(x)$ defined in $\mathbb{R}^n$ is a function such that for every fixed $x \in \mathbb{R}^n$, $T(x)$ is a random variable on the probability space. 
\end{definition}

\begin{definition}
The  covariance function $R(x,y)$ of a random field $T$ is defined as 
$$R(x,y)= \mathbb{E}  \big[ T(x)- \mathbb{E} T(x) \big] \big[ T(y)- \mathbb{E} T(y) \big].$$
\end{definition}

Consider a random field $T$. If  the joint distribution $$F_{x_1,\cdots, x_m}(z_1,\cdots, z_m) = P\big[  T(x_1) \leq z_1, \cdots, T(x_m) \leq z_m  \big]$$
is invariant under the translation $$(x_1,\cdots,x_m) \to (x_1+\tau, \cdots, x_m +\tau),$$ $T$ is said to be stationary or homogeneous. For a stationary random field $T$, we can show 
$$\mathbb{E} T(x) = \mathbb{E} T(0)$$ 
and subsequently
$$R(x,y)=f(x-y)$$ for some function $f$. Although the converse is not always true, such a case is not often encountered in practical applications \citep{yaglom.1987}  so we may equate the stationarity with the condition $$\mathbb{E} T(x) = \mathbb{E} T(0), \; R(x,y)=f(x-y).$$ 
A special case of stationary fields is an isotropic field which requires the covariance function to be rotation invariant, i.e. $$R(x,y) = f(\| x - y\|)$$ for some function $f$. $\| \cdot \|$ is the geodesic distance in the underlying manifold. 

\index{covariance fields}
\index{fields ! covariance}

\subsection{Gaussian Fields}

\index{fields ! Gaussian}
\index{Gaussian fields}

The most important class of random fields is Gaussian fields. A more rigorous treatment can be found in \citet{adler.2007}. Let us start defining a multivariate normal distribution from a Gaussian random variable. 

\begin{definition}
A random vector $T=(T_1,\cdots, T_m)$ is multivariate normal if 
$\sum_{i=1}^m c_i T_i$ is Gaussian for every possible $c_i \in \mathbb{R}$.
\end{definition}
Then a Gaussian random field can be defined from a multivariate normal distribution.
\begin{definition}
\label{def:gaussian}
A random field $T$ is a {\em Gaussian random field} if $T(x_1),\cdots, T(x_m)$ are 
multivariate normal for every $(x_1,\cdots, x_m) \in \mathbb{R}^m$.  
\end{definition}
An equivalent definition to Definition \ref{def:gaussian} is as follows. 
\begin{definition}
$T$ is a  Gaussian random field if the finite joint distribution\\ 
$F_{x_1,\cdots, x_m}(z_1,\cdots, z_m)$ is a multivariate normal for every $(x_1,\cdots, x_m)$. 
\end{definition}

$T$ is a mean zero Gaussian field if $\mathbb{E}T(x) = 0$ for all $x$. Because any mean zero multivariate normal distribution can be completely characterized by its covariance matrix, a mean zero Gaussian random field $T$ can be similarly  determined by its covariance function $R$. Two fields $T$ and $S$ are
independent if $T(x)$ and $S(y)$ are independent for every $x$
and $y$. For mean zero Gaussian fields $T$ and $S$, they are
independent if and only if the cross-covariance function 
$$R(x,y)=\mathbb{E} \big[ T(x)T(y) \big]$$
vanishes for all $x$ and $y$.

Given two arbitrary mean zero
Gaussian fields, is there mapping that makes them independent? Let
$e=(e_1(t),e_2(s))^{\top}$ be a vector field. Let $A$ be a constant
matrix. Consider transformation $Ae$ and its covariance
$$\mathbb{E} [Ae (Ae) ^{\top}] = A\mathbb{E} [ee ^{\top}] A ^{\top}.$$ 
Note that
$\mathbb{E} [ee ^{\top}]$ is symmetric and its diagonal terms are
positive so it is a symmetric positive definite matrix so we have
a singular value decomposition of the form $\mathbb{E}[ee ^{\top}] =
QDiag(\lambda_1,\lambda_2)Q ^{\top}$ where $Q$ is an orthogonal matrix.
Simply let $A=Q ^{\top}$ and it should make the component of $Ae$
uncorrelated for all $t$ and $s$. But they are still not
independent.

The Gaussian white noise is a Gaussian random field with the Dirac-delta function $\delta$ as the covariance function. Note the Dirac delta function is defined as $\delta(x) = \infty, x =0$, $\delta(x) = 0 x \neq 0$ and $\int \delta(x) = 1$. Numerically we can simulate the Dirac delta function as the limit of the sequence of Gaussian kernel $K_{\sigma}$ when $\sigma \to \infty$. The Gaussian white noise is simulated as an independent and identical Gaussian random variable at each voxel.

\subsection{Derivative of Gaussian Fields} 

\index{fields ! derivatives}

Any linear operation $f$ on Gaussian fields is again Gaussian fields.
Suppose $\mathcal{G}$ be a collection of Gaussian random fields. Then $f(\mathcal{G}) \subset \mathcal{G}$. For given $X, Y \in \mathcal{G}$, we have $c_1 X + c_2 Y \in \mathcal{G}$ again for all $c_1$ and $c_2$. Therefore, $\mathcal{G}$ forms an infinite-dimensional vector space. Not only the linear combination of Gaussian fields is again Gaussian but also the derivatives of Gaussian fields are Gaussian. To see this, we define the mean-square convergence. 

\begin{definition}
A sequence of random fields $T_h$, indexed by $h$ converges to
$T$ as $h \to 0$ in mean-square if
$$\lim_{h \to 0} \mathbb{E} \big |T_h - T  \big|^2=0.$$
\end{definition}
We will denote the mean-square convergence using the usual limit notation:
$$\lim_{h \to 0} T_h = T.$$
The convergence in mean-square implies the convergence in mean. This can be seen from 
$$\mathbb{E} \big|T_h  - T \big|^2 =  \mathbb{V} \big[T_h - T \big]^2 +\big( \mathbb{E} | T_h - T | \big)^2.$$
Now let $T_h \to T$ in mean square. Each term in the right hand side should also converges to zero proving the statement.

Now we define the derivative of field in mean square sense as
$$\frac{d T(x)}{d x} = \lim_{h \to 0} \frac{T(x+h)-T(x)}{h}.$$
Note that if $T(x)$ and $T(x+h)$ are Gaussian random fields, $T(x+h)-T(x)$ is again Gaussian, and hence the limit on the right hand side is again Gaussian. If $R$ is the covariance function of the mean zero Gaussian field $T$, the covariance function of its derivative field is given by
$$\mathbb{E} \Big[\frac{d T(x)}{d x}
\frac{d T(y)}{d y} \Big]=\frac{\partial ^2 R(x,y)}{\partial x
\partial y}.$$

Given zero mean Gaussian field $X(t), t =(t_1,\cdots,t_n) \in \mathbb{R}^n$, the Hessian  field $H(t)$ of $X(t)$ is given by 
$$H(t) = \Big(\frac{\partial^2
X(t)}{\partial t_i \partial t_j}\Big).$$ 
If we have mean zero Gaussian random variables $Z_1,\cdots, Z_n, \mathbb{E} (Z_1 \cdots
Z_n) = 0$ if $n$ odd. Hence the expectation of the determinant of
Hessian of a mean zero Gaussian field vanishes. For $n$ even and
assuming isotropic covariance $R(t,s) = R(\|t-s\|)$, i.e.,
stationarity, we can further simplify the expression.

\subsection{Integration of Gaussian Fields}

\index{fields ! integration}

The integration of Gaussian fields is also Gaussian. To see this, define the integration of a random field as the limit of
Riemann sum. Let $\cup_{i=1}^n \mathcal{M}_i$ be a partition of $\mathcal{M}$, i.e. 
$$\mathcal{M} = \cup_{i=1}^n \mathcal{M}_i \mbox{ and } 
\mathcal{M}_i \cap \mathcal{M}_j = \emptyset \mbox{ if } i \neq j.$$ 
Let $x_i \in
\mathcal{M}_i$ and $\mu(\mathcal{M}_i)$ be the volume of $\mathcal{M}_i$. Then we define the integration of field $T$ as
$$\int_{\mathcal{M}} T(x) \; dx = \lim_{} \sum_{i=1}^n T(x_i)\mu(\mathcal{M}_i),$$
where the limit is taken as $\mu(\mathcal{M}_j) \to 0$ for all $j$.

Multiple integration is defined similarly. When we integrate a Gaussian field, it is the limit of a linear combination of
Gaussian random variables so it is again a Gaussian random variable. In general, any linear operation on Gaussian fields will result in  a Gaussian field. 

Let $X(t)$ be a zero mean Gaussian fields in
$\Omega$ with covariance function $R$. Let's find the distribution of
$\int_{\Omega} X(t) \; dt$. 
Obviously this is a zero mean Gaussian random variable so we only need to
find the second moment 
\bq &&\mathbb{E} \big[\int_{\Omega} X(t)
\;dt \int_{\Omega} X(s) \;ds\big]\\
&=& \int_{\Omega}\int_{\Omega} \mathbb{E} [X(t)X(s)] \;dt \;ds\\
&=& \int_{\Omega}\int_{\Omega} R(t,s) \;dt \;ds.
\eq

\subsection{$t$, $F$ and $\chi^2$ Fields} 

\index{fields ! $t$-fields}
\index{fields ! $F$-fields}
\index{fields ! $\chi^2$-fields}

We can use i.i.d. Gaussian fields to construct  $\chi^2$-, $t$-, $F$-fields, all of which are extensively studied in \citep{cao.1999.stat, worsley.2004, worsley.1996, worsley.1994}. The $\chi^2$-field with $m$ degrees of freedom is defined as
$$T(x) = \sum_{i=1}^m X_i^2(x),$$
where $X_1,\cdots,X_m$ are independent, identically distributed Gaussian fields with zero mean and unit variance. Similarly, we can define $t$ and $F$ fields as well as Hotelling's $T^2$ field. The Hotelling's $T^2$-statistic has been widely used in detecting morphological changes in deformation-based morphometry \citep{cao.1999.stat, collins.1998, gaser.1999, joshi.1998, thompson.1997}. In particular,  \citep{cao.1999.stat} derived the excursion probability of the Hotelling's $T^2$-field and applied it to detect gender specific morphological differences.

\section{Convolution on Random Fields}
Consider the following integral
$$Y(t) = \int K(t,s)X(s) \; ds.$$
where $K$ is called the kernel of the integral.  Define {\em convolution} between kernel $K$ and random field $X$ as the above integral.
\bq Y(t) = K*X(t).\label{eq:lecture2-con}\eq 

Suppose  the kernel to be isotropic probability density, i.e. $K(t,s) =
K(t-s)$ and $\int K(t) \; dt =1.$ Further we may assume $K$ to be
unimodal with some parameter $\sigma$ such that
$$ \lim_{\sigma \to 0} K(t;\sigma) \to \delta(t),$$
the {Dirac-delta function}. Since the Dirac-delta function satisfies 
$$\int \delta(t-s)f(s) \;ds = f(t),$$
it can be easily seen that
$$\lim_{\sigma \to 0} K(\cdot;\sigma)*X \to X.$$

\subsection{Kernel smoothing estimator} 
Let $t=(t_1,\cdots,t_n)' \in
\mathbb{R}^n.$ The $n$-dimensional isotropic Gaussian kernel is
given by the products of $1$-dimensional Gaussian kernel:
$$K(t)=\prod_{i=1}^n \frac{1}{\sqrt{2\pi}}e^{-t_i^2/2}.$$
The isotropic kernel under linear transform $t \to Ht$ changes the
shape of the kernel to anisotropic kernel
$$K_H(t) = \frac{1}{det (H)} K(H^{-1}t).$$
Note the multivariate kernel $K_H$ is a probability distribution, i.e.,
$$\int  K_H(t) \; dt = 1.$$
$(\det H)^{-1}$ is the Jacobian of the transformation that normalize
the density. It is the distribution of $n$-dimensional
multivariate normal with the covariance matrix $HH'$, i.e.
$N(0,HH')$. $H$ is also called the {\em bandwidth matrix} in the
context of kernel smoothing and it measures the amount of
smoothing. It can be shown that $K_{H}$ satisfy the definition of
the Dirac-delta function as the eigenvalues $\lambda_i$ of $H$ go
to zero, i.e.
$$\lim_{\lambda_i \to 0} K_{H}(t) = \delta(t).$$
The limit of the sequence of the Gaussian kernels gives the
Dirac-delta function and this is how we implement the Dirac-delta
function in computer programs. From now on we let $H \to 0$ if all
$\lambda_i \to 0$ and $H \to \infty$ if all $\lambda_i \to
\infty$.

Suppose we have an additive model \bq  Y(t) = \mu(t) + \epsilon(t) \label{eq:model1}\eq where
$e$ is a mean zero random field and $\mu$ is an unknown signal. In
image analysis, observations are so dense that we can take them to
be continues functional data. Then the kernel smoothing estimator
is given by 
\bqn \widehat \mu(t) = K_{H}*Y(t)= \int K_{H}(t-s) Y(s)\;ds\label{eq:model2}\eqn  
As $H$ goes to zero, we are smoothing
less. To see this note that
$$\lim_{\sigma \to 0} \widehat \mu(t) = \int \delta(t-s) Y(s) \; ds = Y(t)$$
 From (\ref{eq:model2}),
$$\mathbb{E} \widehat \mu(t) = K_{H}*\mu(t) \to \mu(t) \mbox { as } \sigma \to 0$$
From the property of the Dirac-delta, as $H \to 0$, $K_{H}*\mu(t)$
converges to the true but unknown signal $\mu$. So
we can see that our kernel estimator becomes more unbiased as $H
\to 0$.

We can show that the estimator $ \widehat \mu(t)$ is a
solution to a heat equation and the condition $H \to \infty$ is
equivalent to the steady state that is reached when we diffuse heat for infinite
amount of time \citep{chung.2012.CNA}.

Other properties of kernel smoothing estimator is as follows.
Assuming $|\mu| \leq \infty$, \bqn \mathbb{E} \widehat \mu (t) \leq
\int K_H(t) \sup_{} \mu(t) \;dt \leq \sup_{} \mu(t).\eqn

Similarly we can bound from below so that 
\bqn \inf_{} \mu(t) \leq
\mathbb{E} \widehat \mu(t) \leq \sup_{} \mu(t) .
\label{eq:lecture3-inf}\eqn

 Inequality (\ref{eq:lecture3-inf}) implies that smoothed signal will be smaller than the maximum and larger than
the minimum of the signal in average. Other interesting property is 
\bqn \int K_H*Y(t) \;
dt &=& \int_t \int_s K_H(t-s)Y(s) \; ds\; dt\\
&=& \int Y(s) \; ds.\eqn
Thus we have
$$\mathbb{E} \int \widehat \mu(t) \; dt = \int \mu(t) \;dt.$$ 
If $\mu(t)$ is a probability such as the transition probability in
Brownian motion, it implies that the total probability is
conserved even after smoothing.

Let $R_Y$ be the covariance function of field $Y$ with $\mathbb{E} Y =0$.
It is trivially $$R_Y(t,s) =
R_{\epsilon}(t,s) =\mathbb{E} [\epsilon(t)\epsilon(s)],$$
the covariance function of $\epsilon$. Then we can show that the covariance function of the
kernel smoother can be shown to be \bqn R(t,s) =\int \int
K_{H}(t-t')K_{H}(t-s')R_{\epsilon}(t',s') \;dt'
 \;ds'.\eqn

\section{Numerical Simulation of Gaussian Fields}
For the random field theory based multiple correction to work, it is necessary to have smooth images. In this section, we show how to simulate smooth Gaussian fields by performing Gaussian kernel smoothing on white noise. This is the easiest way of simulating Gaussian fields.

White noise is defined as a random field whose covariance function is proportional to the Dirac-delta function $\delta$, i.e. 
$$R(x,y) \propto \delta(x-y).$$
For instance, me may take $R(x,y) = \lim_{\sigma \to 0} K_{\sigma}(\|x
- y\|)$, the limit of the usual isotropic Gaussian kernel. White noise is usually characterized via generalized functions.

One example of white noise is the generalized derivative of
Brownian motion (Wiener process) called white Gaussian noise. 

\begin{definition}
\label{def:brownian}
Brownian motion (Wiener process) $B(x), x \in \mathbb{R}^+$ is zero mean Gaussian
field with covariance function $$R_B(x,y) = \min(x,y).$$ 
\end{definition}

Following Definition \ref{def:brownian}, we can show ${\bf Var} B(x) = x$ by taking $x=y$ in the covariance function and $B(0)=0$ by letting $x=0$ in the variance. The increments of Wiener processes in nonoverlapping
intervals are independent identically distributed (iid) Gaussian. Further the paths of the Wiener process
is continous in the Kolmogorov sense while it is not differentiable. For a different but identical canonical
construction of Brownian motion, see \citep{oksendal.2010}. 
Higher dimensional
Brownian motion can be generalized by taking each component of
vector fields to be i.i.d. Brownian motion.

\index{Brownian motion}
\index{Gaussian noise}

\begin{figure}[t]
\centering
\includegraphics[width=0.9\linewidth]{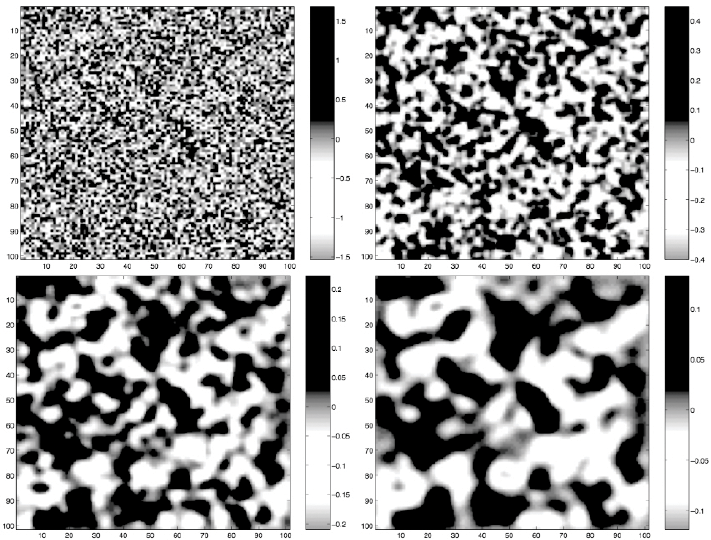}
\caption{Random fields simulation via iterated Gaussian kernel
smoothing with $\sigma=0.4$. $N(0,0.4^2)$. White noise, 1, 4 and 9
iterations in sequence.\label{fig:lec9-1}}
\end{figure}

Although the path of Wiener process is not differentiable, we can define the generalized derivative via integration by parts with a smooth function $f$
called a test function in the following way
$$f(x)B(x) = \int_0^x f(y)\frac{dB(y)}{dy} \;dy + \int_0^x \frac{f(y)}{dy}B(y)\;dy.$$
Taking the expectation on both sides we have
$$ \int_0^x f(y) \mathbb{E}\frac{dB(y)}{dy} \;dy =0.$$
It should be true for all smooth $f$ so $\mathbb{E}
\frac{dB(y)}{dy} =0$. Further it can be shown that the covariance
function of process $dB(y)/dy \propto \delta(x-y)$.

\index{Wiener process}

The Gaussian white noise can be used to construct smooth Gaussian
random fields of the form
$$X(x) = K * W(x) = K * \frac{dB(x)}{dx},$$
where $K$ is a Gaussian kernel. Since Brownian motion is zero mean Gaussian process, $X(x)$ is
obviously zero mean field with the covariance function 
\bqn
R_X(x,y) &=& \mathbb{E} [K * W(x)K * W(y)]\\
&\propto&\int K(x-z)K(y-z) \; dz.\eqn The case when $K$  is an 
isotropic Gaussian kernel was investigated by D.O.
Siegmund and K.J. Worsley with respect to optimal filtering in
scale space theory \citep{siegmund.1996}.

In numerical implementation, we use the discrete white Gaussian noise which is simply a Gaussian random variable. 

\begin{example}
Let $w$ be a discrete version of white Gaussian noise given by
$$w(x) = \sum_{i=1}^m Z_i \delta(x-x_i),$$
where i.i.d. $Z_i \sim N(0,\sigma_w^2)$. Note that
\bqn K* w(x) = \sum_{i=1}^m Z_iK(x-x_i). 
\label{eq:fields-KW}\eqn
The collection of random variables $K* w(y_1), \cdots, K*w(y_l)$ forms a
multivariate normal at arbitrary points $y_1,\cdots, y_l$. Hence
the field $K * w(x)$ is a Gaussian field. 
\end{example}
The covariance function of the field \ref{eq:fields-KW} is given by
\bqn R(x,y) &=& \sum_{i,j=1}^m \mathbb{E}(Z_iZ_j) K(x-x_i)
K(y-x_j) \\
&=&\sum_{i=1}^m \sigma_w^2K(x-x_i) K(y-x_i).
 \eqn
 As usual we may take $K$ to be a
Gaussian kernel. Let us simulate some Gaussian fields. 

\begin{example}
The unknown signal is assumed to be $\mu(x,y) = \cos (10x) + \sin(8y), (x,y) \in [0,1]^2$ and
white noise error $w \sim N(0,0.4^2)$ which is shown in the top left of Figure \ref{fig:lec9-1}. Then iteratively more smooth version of Gaussian random fields are constructed by

\begin{verbatim}
w=normrnd(0,0.4,101,101);
smooth_w=w;
for i=1:10
  smooth_w=conv2(smooth_w,K,'same');
  figure;imagesc(smooth_w);colorbar;
end;
\end{verbatim}
with kernel weight {\tt K}.
\end{example}

\section{Statistical Inference on Fields}
\label{sec:MCC}
\index{fields ! statistical inference}

Given functional measurement $Y$, we have model
$$Y(x) =\mu(x) + \epsilon(x),$$
where $\mu$ is unknown signal and $\epsilon$ is a zero mean unit variance Gaussian field. We assume $x \in \mathcal{M} \subset \mathbb{R}^n$. In brain imaging, one of the most important problems is that of signal detection, which can be stated as the problem of identifying the regions of
statistical significance. So it can be formulated as an
inference problem $$H_0: \mu(x)=0 \mbox{ for all } x \in \mathcal{M} \mbox{ vs. }\;
H_1:\mu(x)
> 0 \mbox { for some } x \in \mathcal{M}.$$
Let
$$H_0(x):\mu(x) = 0$$ 
at a fixed point $x$. Then the null hypothesis $H_0$ is a collection of multiple hypotheses $H_0(x)$ over all $x$. Therefore, we have
$$H_0 = \bigcap_{x \in \mathcal{M}} H_0(x).$$
We may assume that $\mathcal{M}$ is the region of interest consisting of the finite number of voxels. We also have the corresponding point-wise alternate hypothesis $$H_1(x): \mu(x)
>0$$
and the alternate hypothesis $H_1$ is constructed as
$$H_1 = \bigcup_{x \in \mathcal{M}} H_0(x).$$
If we use $Z$-statistic as a test statistic, for instance, we will reject
each $H_0(x)$ if $Z > h$ for some threshold $h$.  So at each fixed $x$, for level $\alpha
= 0.05$ test, we need to have $h=1.64$.
However, if we threshold at $\alpha =0.05$, $5\%$ of observations are false positives. 
Note that the false positives are pixels where we are incorrectly rejecting $H_0(x)$ when it is actually true. 
However, these are the false positives related to testing $H_0(x)$. For determining the true false positives associated with testing $H_0$, we need to account for multiple comparisons. 

\begin{definition}

\index{type-I error}
\index{error ! type-I}

The type-I error is the probability of  rejecting the null hypothesis (there is no signal) when the alternate hypothesis (there is a signal) is true. 
\end{definition}
The type-I error is also called the {\em family-wise error rate} (FWER) and given by
\bqn \alpha &=& P(\mbox{
reject } H_0 \; | \; H_0 \mbox{ true }) \nonumber \\
&=& P(\mbox{ reject some } H_0(x) \;| \; H_0 \mbox{ true }) \nonumber\\
&=& P\Big(\bigcup_{x \in \mathcal{M}} \{Y(x) > h\} \; \Big| \; \mathbb{E} Y=0\Big). \label{eq:deform-alphacup}
\eqn
Unfortunately, $Y(x)$ is correlated over $x$ and it makes the
computation of type-I error almost intractable for random fields
other than Gaussian.\\

\subsection{Bonferroni Correction} 

\index{Bonferroni correction}
\index{multiple comparisons ! Bonferroni}

One standard method for dealing
with multiple comparisons is to use the Bonferroni correction. Note
that the probability measure is additive so that for any event $E_j$, we have \bq
P \Big(\bigcup_{j=1}^{\infty} E_j \Big) \leq \sum_{j=1}^{\infty} P(E_j).\eq This
inequality is called Bonferroni inequality and it has been used in the construction of simultaneous
confidence intervals and multiple comparisons when the number of
hypotheses are small. From (\ref{eq:deform-alphacup}), we have 
\bqn \alpha &=& P\Big(\bigcup_{x \in \mathcal{M}} \{Y(x) >
h\} \; \Big| \; \mathbb{E} Y=0\Big) \label{eq:deform-bonferoni1}\\
&\leq& \sum_{x \in \mathcal{M}} P\big(Y(x_j)>h \; | \; \mathbb{E} Y
=0\big)\label{eq:deform-bonferoni2}\eqn So by controlling each type-I error
separately at
$$P\big(Y(x_j)>h \; | \; \mathbb{E} Y =0\big) < \frac{\alpha}{\# \mathcal{M}}$$
we can construct the correct level $\alpha$ test. Here $\# \mathcal{M}$ is the number of voxels.

The problem with the Bonferroni correction is that  it is too
conservative. The Bonferroni inequality  (\ref{eq:deform-bonferoni2}) becomes exact when the measurements across voxels are all independent, which is unrealistic. Since the measurements are expected to be  strongly correlated across voxels, we have highly correlated statistics. So in a sense, we have a less number of comparisons to make. 

\begin{figure}[t]
\includegraphics[width=1\linewidth]{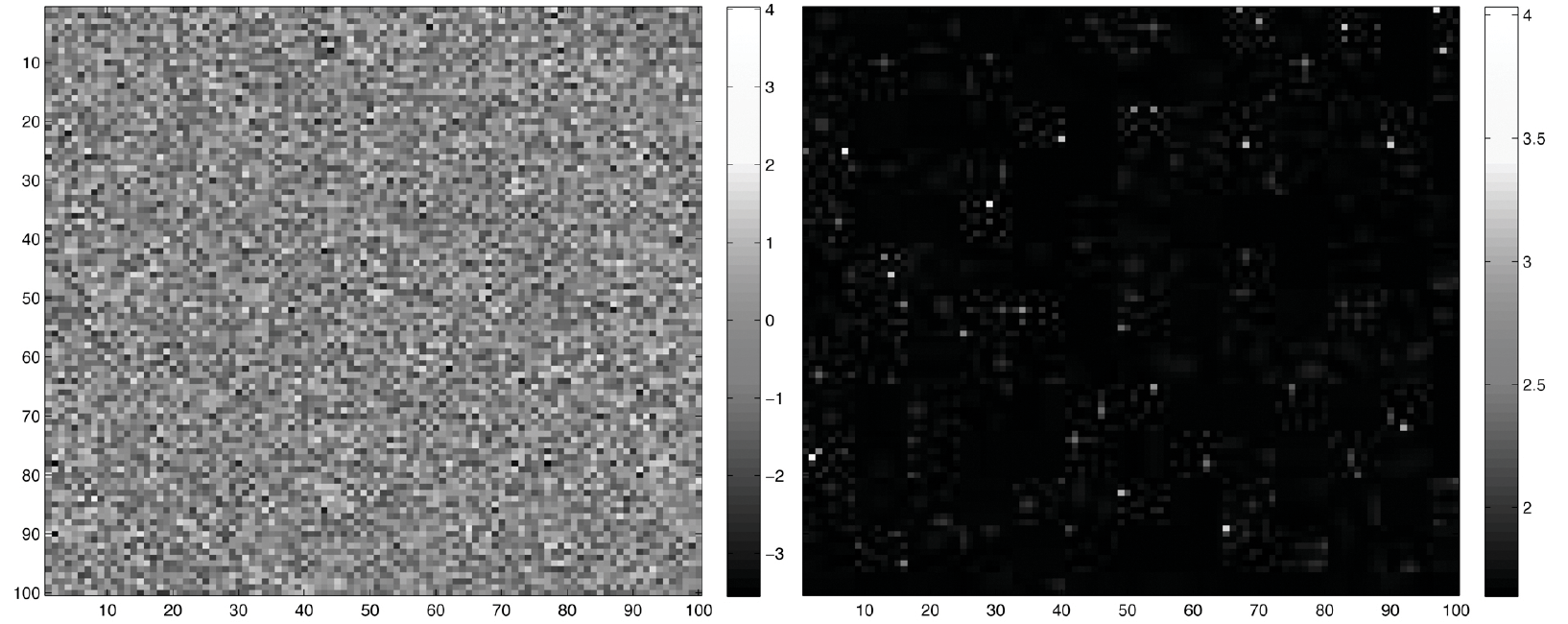}
\centering
\caption{Image consisting of $N(0,1)$ noise. At the thresholding 1.64 corresponding to the significance level $\alpha = 0.05$, 5$\%$ of all pixels are false positives.}
\label{fig:bonferroni}
\end{figure}

Let us illustrate the Bonferroni correction procedure using {\tt MATLAB}. Consider 
$100 \times 100$ image $Y$ of standard normal random variables (Figure \ref{fig:bonferroni}). The threshold corresponding to the significance $\alpha = 0.05$ is 1.64. 

\begin{verbatim}
Y=normrnd(0,1,100,100); 
figure; imagesc(Y); colorbar; colormap('hot')

norminv(0.95,0,1)

ans =

    1.6449
    
[Yl, Yh] = threshold_image(Y, 1.64);
figure; imagesc(Yh); colormap('hot'); colorbar;     
\end{verbatim}

By thresholding the image at 1.64, we obtain approximately about 5$\%$ of pixels as false positives. To account for the false positives, we perform the Bonferroni correction. For image of size $100 \times 100$, there are $10000$ pixels. Therefore, $\alpha/\# \mathcal{M} =0.05/ 10000= 0.000005 $ is the corresponding point-wise $p$-value and the corresponding threshold is 4.42. In this example,  there is no pixel that is higher than 4.42 so we are not detecting any false positives as expected. 

\begin{verbatim}
n=100*100;
size(find(reshape(Y,n,1)>=1.64),1)/n
norminv(1-0.05/10000,0,1)

ans =

    4.4172
\end{verbatim}

\subsection{Rice Formula} 
\index{Rice formula}
\index{multiple comparisons ! Rice formula}

We can obtain a less conservative estimate for (\ref{eq:deform-alphacup}) using the random field theory. Assuming $\mathbb{E} Y=0$, we have
\bqn \alpha(h) &=& P\Big(\bigcup_{x \in \mathcal{M}} \{Y(x) > h\}\Big) \nonumber \\
&=&1-  P\Big(\bigcap_{x \in \mathcal{M}} \{Y(x) \leq h  \}\Big)  \nonumber \\
&=& 1 - P \Big(\sup_{x \in \mathcal{M}} Y(x) \leq h \Big)  \nonumber \\
&=& P \Big(\sup_{x \in \mathcal{M}} Y(x) > h \Big). \label{eq:deform-supYx} \eqn 
In order to construct the $\alpha$-level test corresponding to $H_0$, we need to know the distribution of
the supremum of the field $Y$.  
The corresponding $p$-value based on the supremum of the field, i.e. $\sup_{x \in \mathcal{M}} Y$, is called the {\em corrected $p$-value} to distinguish it from the usual $p$-value obtained from the statistic $Y$. Note that the p-value is the smallest $\alpha$-level at which the null hypothesis $H_0$ is rejected.

Analytically computing the exact distribution of the supremum of random fields is hard. If we denote $Z=\sup_{x \in \mathcal{M}} Y(x)$ and $F_Z$ to be the cumulative
distribution of $Z$, for the given $\alpha =0.05$, we can compuate
$h=1-F_Z^{-1}(\alpha)$. Then the region of statistically
significant signal is localized as $\{x \in \mathcal{M}: Y(x) > h\}$.

The distribution of supremum of Brownian motion is somewhat simple due to its independent increment properties. However, for smooth random fields, it is not so straightforward. Read \citep{adler.2000} for an overview of computing the distribution of the supremum of smooth fields. 

Consider 1D smooth stationary Gaussian random process $Y(x), x \in \mathcal{M}=[0,1] \subset
\mathbb{R}$. Let $N_h$ to be the number of times $Y$ crosses over
$h$ from below (called upcrossing) in $[0,1]$. Then we have \bq
P\Big(\sup_{x \in
[0,1]} Y(x) > h\Big) &=& P(N_h \geq 1 \mbox{ or } Y(0) > h)\\
&\leq& P(N_h \geq 1) + P(Y(0) >h)\\
&\leq& \mathbb{E} N_h + P(Y(0) >h).\eq
If $R$ is the covariance function of the field $Y$, we have 
$$R(0) = \sigma^2 = \mathbb{E} Y^2(x).$$
It can be shown that from Rice formula \citep{adler.1993, rice.1944},
$$\mathbb{E} N_h = \frac{1}{\pi}\left(   \frac{-R''(0)}{R(0)} \right)^{1/2} \exp \left(\frac{h^2}{2\sigma^2}\right).$$
Also note that 
$$P(Y(0) > h) =1-\Phi \Big(\frac{h}{\sigma} \Big)$$ where $\Phi$ is the
cumulative distribution function of the standard normal. Then from the inequality that bounds the cumulative distribution of the standard normal \citep{feller.1968}, we have \bq
\Big( 1-\frac{\sigma^2}{h^2} \Big)\frac{\sigma}{\sqrt{2\pi}
h}e^{-h^2/2\sigma^2} \leq 1-\Phi \Big(\frac{h}{\sigma} \Big) \leq
\frac{\sigma}{\sqrt{2\pi} h}e^{-h^2/2\sigma^2} \eq
So \bq P\Big(\sup_{x \in [0,1]} Y(x) > h\Big) \leq \Big[c_1
+\frac{c_2}{\sqrt{2\pi} h}\Big] e^{-h^2/2\sigma^2}\eq
for some $c_1$ and $c_2$.  In fact we can show
that 
\bq P\Big(\sup_{x \in [0,1]} Y(x) > h\Big) = \Big[c_1
+\frac{c_2}{h} +O(h^{-2})\Big] e^{-h^2/2\sigma^2}.\eq
\\

\subsection{Poisson Clumping Heuristic}  

\index{multiple comparisons ! Poisson clumping heuristic}

To extend the Rice formula to a higher dimension, we need a different mathematical machinery.  For this method to work, the random field $Y$ needs to be sufficiently smooth and isotropic. The smoothness of a random field corresponds to the random field being differentiable. There are very few cases for which exact formulas for the excursion probability (\ref{eq:deform-supYx}) are known \citep{adler.1990}. For this reason, approximating the excursion probability is necessary for most cases.

From the {\it Poisson clumping heuristic} \citep{aldous.1989}, 
$$P\Big( \sup_{x \in \mathcal{M}} Y(x) <  h \Big) \approx \exp \biggl(-\frac{\| \mathcal{M} \|}{\mathbb{E} \|A_h\|} P\big( Y(x) \geq h\big)\biggl),$$
where $\| \cdot \|$ is the Lebesgue measure of a set and the random set 
$$A_h = \{ x \in \mathcal{M}
: Y(x) > h\}$$
is called the {\em excursion set} above the threshold $h$. This approximation involves unknown $\mathbb{E} \|A_h\|$, which is the mean clump size of the excursion set. The distribution of $\|A_{h}\|$ has been estimated for the case of Gaussian \citep{aldous.1989}, $\chi^2 ,t$ and $F$ fields \citep{cao1999.AAP} but for general random fields, no approximation is available yet.



\section{Expected Euler Characteristics} 

\index{multiple comparisons ! Euler characteristic method}
\index{Euler characteristics}

An alternate approximation to the supremum distribution based on the expected Euler characteristic (EC) of $A_h$ is also available. The Euler characteristic approach reformulates the geometric problem  as a topological problem. Read \citep{adler.1981, cao.2001, taylor.2007, worsley.2003} for an overview of the Euler characteristic method. 

\begin{figure}[t]
\includegraphics[width=1\linewidth]{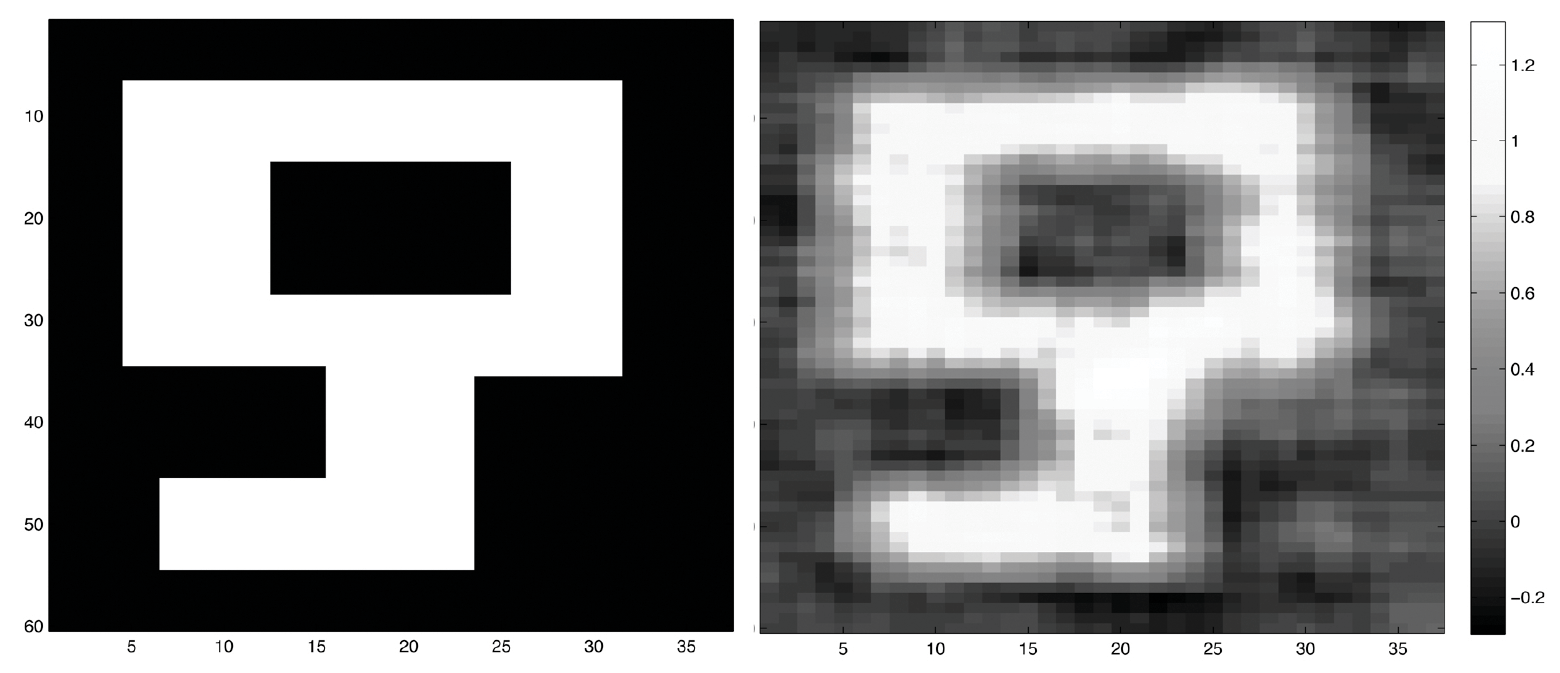}
\centering
\caption{Gaussian white noise is added and smoothed to the key shaped binary object. The Euler characteristic of an object with a hole is 0.}
\label{fig:power-euler}
\end{figure}

\begin{figure}[t]
\includegraphics[width=0.8\linewidth]{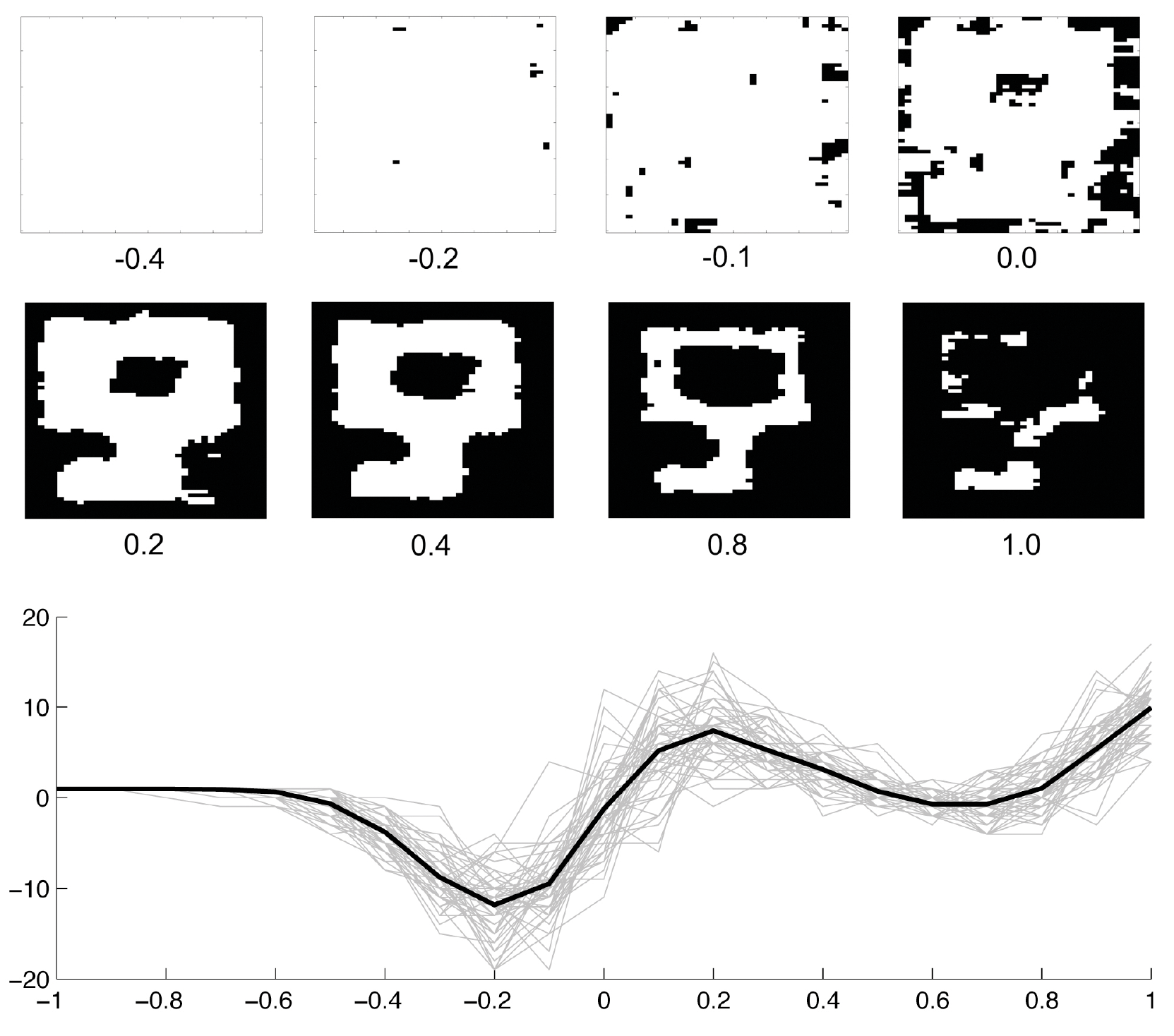}
\centering
\caption{The mean Euler characteristic of the excursion sets of the shaped object shown in Figure \ref{fig:power-euler}. The number below each object is the threshold. \label{fig:power-eulerplot}}
\end{figure}

For sufficiently high threshold $h$, it is known that
\bqn P \Big(\sup_{x \in \mathcal{M}} Y(x) > h \Big) \approx \mathbb{E} \chi(A_h) = \sum_{d=0}^N
\mu_d(\mathcal{M})\rho_d(h) 
\label{eq:deform-ECexpansion}
\eqn 
where $\mu_d(\mathcal{M})$ is the $d$-th
Minkowski functional or {\em intrinsic volume} of $\mathcal{M}$ and
$\rho_d$ is the $d$-th Euler characteristic (EC) density of
$Y$ \citep{worsley.1998}. For details on intrinsic volume, read \citep{schmidt.2005}. The expansion (\ref{eq:deform-ECexpansion}) also holds for non-isotropic fields but we will not pursue it any further. 
Compared to other approximation methods such as the Poisson clump heuristic and the tube formulae, the advantage of using the Euler characteristic formulation is that a simple exact expression can be found for $\mathbb{E}\; \chi(A_h)$. Figure \ref{fig:power-euler} and Figure \ref{fig:power-eulerplot} show how $\chi(A_h)$ and $\mathbb{E}\; \chi(A_h)$ change as the threshold $h$ increases for a simple binary object with a hole.

\subsection{Intrinsic Volumes}  
\index{intrinsic volumes}

The $d$-th intrinsic volume of $\mathcal{M}$ is a generalization of $d$-dimensional volume. Note that $\mu_0(\mathcal{M})$ is the Euler characteristic of $\mathcal{M}$. 
$\mu_N(\mathcal{M})$ is the volume of $\mathcal{M}$ while $\mu_{N-1}(\mathcal{M})$ is half the surface area of $\mathcal{M}$. There are various techniques for computing the intrinsic volume \citep{taylor.2007}. The methods depend on the smoothness of the underlying manifold $\mathcal{M}$. For a solid sphere with
radius $r$, the intrinsic volumes are 
$$\mu_0 =1, \mu_1 = 4r, \mu_2= 2 \pi r^2, \mu_3 = \frac{4}{3}\pi r^3.$$ For a 3D box of size $a \times b \times c$, the intrinsic volumes are 
$$\mu_0=1, \mu_1= a+ b+c, \mu_2 = ab + bc + ac, \mu_3=abc.$$ In general, the intrinsic volume can be given in terms of a curvature matrix. Let $K_{\partial \mathcal{M}}$ be the curvature matrix of $\partial \mathcal{M}$ and $\text{detr}_{d}(K_{\partial \mathcal{M}}) $ be the sum of the determinant of all  $d \times d$ principal minors of $K_{\partial \mathcal{M}}$. For $d=0,\cdots,N-1$ the Minkowski functional $\mu_d(\mathcal{M})$ is defined as 
$$\mu_d(\mathcal{M})=\frac{\Gamma(\frac{N-i}{2})}{2\pi^{\frac{N-i}{2}}}\int_{\partial \mathcal{M}} \text{detr}_{N-1-d}(K_{\partial \mathcal{M}}) \; dA,$$
and $\mu_N(\mathcal{M})=\|\mathcal{M}\|$, the Lebesgue measure of $\mathcal{M}$.

For irregular jagged shapes such as the 2D corpus callosum shape $\mathcal{M}$,
the intrinsic volume can be estimated in the following
fashion \citep{worsley.1996.hbm, chung.2004.ni}. Treating pixels inside $\mathcal{M}$ as points on a lattice, let $V$ be the
number of vertices that forms the corners of pixels, $E$ be the number of edges connecting each adjacent lattice points and $F$ be the number of faces formed by four
connected edges. We assume the distance between the adjacent lattice points is $\delta$ in all directions. Then $$\mu_0=V - E + F, \mu_1=(E-2F)\delta, \mu_2= F\delta^2.$$  To find the number of edges and pixels
contained in $\mathcal{M}$, we start from an initial face (pixel) somewhere in  the corpus callosum and add one face at a time  while counting the additional edges and faces. In this fashion, we can grow a graph that will eventually contains all the pixels that form the corpus callosum. A numerical method for computing the intrinsic volume 
for jagged irregular shapes has been implemented in FMRISTAT package (\url{www.math.mcgill.ca/keith/fmristat}).


\subsection{Euler Characteristic Density} 
\index{Euler characteristics} 

The $d$-th EC-density is given by
$$\rho_d(h) = \mathbb{E} \big[(Y> h)\det( -\ddot{Y_d})|\dot{Y_d} = 0 \big]
P(\dot{Y_d}=0),$$ where dot notation indicates partial differentiation
with respect to the first $d$ components. The subscript $d$ represents the first $d$ components of $Y$. Computation of the
conditional expectation is nontrivial other than for Gaussian fields.
For zero mean and unit variance Gaussian field $Y$, we have for instance
\bq \rho_0 &=& P(Y > h) = 1 - \Phi(h)\\
\rho_1 &=& \lambda^{1/2}\frac{e^{-h^2/2}}{2\pi}\\
\rho_2 &=& \lambda h\frac{e^{-h^2/2}}{(2\pi)^{3/2}}\\
\rho_3 &=& \lambda^{3/2}(h^2-1) \frac{e^{-h^2/2}}{(2\pi)^2}, \eq
where $\lambda$ measures the smoothness
of fields, defined as the variance of the derivative of component
of $Y$. The exact expression for the EC density $\rho_d$ is available for other random fields such as $t,\chi^2,F$ fields \citep{worsley.1994}, Hotelling's $T^2$ fields \citep{cao.1999.stat} and scale-space random fields \citep{siegmund.1996}. In each case, the EC density $\rho_d$ is proportional to $\lambda^{\frac{d}{2}}$ and it changes depending on the smoothness of the field. 

If $X_1,\cdots,X_{\alpha}, Y_1,\cdots,
Y_{\beta}$ are i.i.d. stationary zero mean unit variance Gaussian fields. Then $F$-field with $\alpha$ and $\beta$ degrees of freedom is given by
$$F(x) = \frac{\sum_{j=1}^{\alpha} X_j^2(x)/\alpha}{\sum_{j=1}^{\beta} Y_j^2(x)/\beta}.$$
To avoid singularity, we need to assume the total degrees of
freedom $\alpha + \beta \gg N$ to be sufficiently larger than the
dimension of space \citep{worsley.1994}. The EC-density
for $F$-field is then given by \bq \rho_0 &=&\int_h^\infty {{\Gamma
({{\alpha + \beta}\over 2})}\over {\Gamma ({\alpha \over 2})\Gamma
({\beta \over 2})}} {\alpha\over \beta}\left({{\alpha x}\over
\beta }\right)^{(\alpha-2) \over 2} \left(1+{{\alpha x}\over \beta
}\right)^{-{(\alpha+\beta)\over 2}} d x,\\
\rho_1&=&\lambda^{1/2} {{\Gamma ({{\alpha+\beta-1}\over
2})2^{1\over 2}}\over{ \Gamma ({\alpha\over 2})\Gamma ({\beta
\over 2})}} \left({{\alpha h}\over \beta}\right)^{(\alpha-1) \over
2} \left(1+{{\alpha h}\over \beta}\right)^{-{(\alpha+\beta-2)\over
2}},\\
\rho_2&=&\frac{\lambda}{2\pi}{{\Gamma ({{\alpha +\beta-2}\over 2})}\over
{\Gamma ({\alpha \over 2})\Gamma ({\beta \over 2})}}
\left({{\alpha h}\over \beta }\right)^{(\alpha-2) \over 2}
\left(1+{{\alpha h}\over \beta }\right)^{-{(\alpha +\beta-2)\over 2}}\\
&\times&\left[(\beta -1){{\alpha h}\over \beta}-(\alpha-1)\right].\\
\eq

If the random field $Y$ is given as the convolution of a smooth kernel $K_h(x)=K(x/h)/h^N$ with a white Gaussian noise \citep{siegmund.1996, worsley.1992}, the covariance matrix of $\dot Y= dY / dx$ is given by
$${\bf Var}(\dot Y) = \frac{\int_{\mathbb{R}^N}\dot K(\frac{x}{h}) \dot K^t(\frac{x}{h}) \;dx}{h^2\int_{\mathbb{R}^N} K^2(\frac{x}{h}) \; dx}.$$
Applying it to a Gaussian kernel $K(x)= (2\pi)^{-n/2}e^{-\|x\|^2/2}$ gives 
$$ \lambda ={\bf Var}(\dot Y_1) = 1/(2h^2).$$ In terms of FWHM of the kernel $K_h$, 
$$\lambda = 4\ln 2 /\text{FWHM}^2.$$

\subsection{Numerical Implementation of Euler Characteristics}
In this section, we show how to compute the expected Euler characteristic in {\tt MATLAB}. The presented routine can be used in estimating the excursion probability numerically. Consider a 2D binary object (Figure \ref{fig:power-euler}), which is stored as a 2D image {\tt toy-key.tif}. After loading the image using {\tt imread}, we perform scaling on image intensity values so that it becomes a binary object. The Euler characteristic of the binary object is then computed using {\tt bweuler}. 

\begin{verbatim}
I=imread('toy-key.tif');
I= imresize(I,.1, 'nearest');
I=(max(max(I))-I);
I=I/max(max(I));
I=double(I);
figure;imagesc(I); colormap('hot')
eul = bweuler(I)

eul =

     0
\end{verbatim}
Since there is a hole in the object, the Euler characteristic is 0. We will add Gaussian white noise $N(0,0.5^2)$ to the binary object and smooth out with FWHM of 10 using {\tt gaussblur}. The resulting image is a Gaussian random field with sufficient smoothness. The smoothed image is stored as {\tt smooth} and displayed in Figure \ref{fig:power-euler}.
\begin{verbatim}
e=normrnd(0,1, 60, 37); 
Y=I + e;
figure;imagesc(Y); colorbar; colormap('hot')
smooth = gaussblur(Y,10);
figure;imagesc(smooth);colorbar; colormap('hot');
\end{verbatim}
At each threshold {\tt h} between -1 and 1, we threshold {\tt smooth} and store it as a new variable {\tt excursion}. Then compute the Euler characteristic of {\tt excursion}. For computing the mean of the Euler characteristic, we simulated Gaussian random fields 50 times using the {\tt for}-loop. 

\begin{verbatim}
figure;
eulsum=zeros(1,21);
for k=1:50
    Y=I+ normrnd(0,1,60,37);
    smooth = gaussblur(Y,10);
    eul=[];
    j=1;
    for h=-1:0.1:1
        [Yl, Yh] = threshold_image(smooth, h);
        Yh=reshape(Yh,60*37,1);
        excursion=zeros(60*37,1);
        excursion(find(Yh>h))=1;
        excursion=reshape(excursion, 60, 37);
        eul(j) = bweuler(excursion);
        j=j+1;
    end;
    hold on; plot(-1:0.1:1, eul, 'Color', [0.7 0.7 0.7])
    eulsum=eulsum+eul;
end;
hold on; plot(-1:0.1:1, eulsum/50, 'Color', 'k', 'LineWidth',2)
\end{verbatim}

\bibliographystyle{agsm} 
\bibliography{reference.2020.07.15}

\end{document}